\documentclass[12pt]{article}

\title{Four consecutive primitive elements in a finite field}

\author{Tamiru Jarso\thanks{The first author is supported by the Australian Defence Science and Technology Group.}\\
Defence Science and Technology Group, Canberra, Australia\\ 
tamiru.jarso@dst.defence.gov.au \\ and \\
{Tim Trudgian\thanks{The second author is supported by the Australian Research Council Future Fellowship FT160100094.}
}\\ 
 School of Science, University of New South Wales, \\
	Canberra, ACT 2610, Australia\\
t.trudgian@adfa.edu.au}

\usepackage{amsmath,amssymb,amsbsy,amsfonts,latexsym,amsopn,amstext,cite,amsxtra,euscript,amscd,bm}
\usepackage{epsfig}
\usepackage{fullpage}

\usepackage[linesnumbered,ruled,vlined]{algorithm2e}
\usepackage{algpseudocode}
\usepackage[english]{babel}
\usepackage[utf8]{inputenc}
\usepackage{mathtools}
\usepackage{float}
\usepackage{algcompatible}
\RequirePackage{mathrsfs} \let\mathcal\mathscr

\usepackage{url}
\usepackage{enumerate}
\usepackage{amsthm}
\usepackage{amsmath}
\usepackage{comment}
\usepackage{amssymb}
\usepackage{booktabs}
\usepackage[colorlinks,linkcolor=blue,anchorcolor=blue,citecolor=blue,backref=page]{hyperref}
\usepackage{amsmath}
\usepackage{amsfonts}
\usepackage{graphicx}
\usepackage[colorinlistoftodos]{todonotes}
\usepackage{pifont}
\usepackage{xcolor}
\usepackage{multirow}

\usepackage[section]{placeins}
\usepackage{float}
\restylefloat{table}

\usepackage{array}    
\usepackage{booktabs} 
\usepackage{tabularx} 
\usepackage{varioref}
\usepackage{afterpage}
\usepackage{seqsplit}

\SetKwProg{Fn}{Function}{}{}
\SetAlFnt{\sffamily}

\SetAlCapFnt{\normalfont\sffamily\large}

\makeatletter
\renewcommand{\algocf@Vline}[1]{
  \strut\par\nointerlineskip
  \algocf@push{\skiprule}
  \hbox{\bgroup\color{cyan}\vrule\egroup%
    \vtop{\algocf@push{\skiptext}
      \vtop{\algocf@addskiptotal #1}\bgroup\color{cyan}\Hlne\egroup}}\vskip\skiphlne
  \algocf@pop{\skiprule}
  \nointerlineskip}
\renewcommand{\algocf@Vsline}[1]{
  \strut\par\nointerlineskip
  \algocf@bblockcode%
  \algocf@push{\skiprule}
  \hbox{\bgroup\color{cyan}\vrule\egroup
    \vtop{\algocf@push{\skiptext}
      \vtop{\algocf@addskiptotal #1}}}
  \algocf@pop{\skiprule}
  \algocf@eblockcode%
}
\makeatother

\newtheorem{theorem}{Theorem}
\newtheorem{lemma}[theorem]{Lemma}

\numberwithin{equation}{section}
\numberwithin{theorem}{section}
\numberwithin{table}{section}

\def\cH{{\mathcal H}}

\def\F{{\mathbb F}}

\def\Z{{\mathbb Z}}

\def\\{\cr}
\def\[{\left[}
\def\]{\right]}
\def\<{\langle}
\def\>{\rangle}

\def\le{\leqslant}
\def\ge{\geqslant}

\def\Hb1{\overline{\cH}_{m}}
\def\Ht1{\widetilde{\cH}_{m}}

\makeatletter
\def\@bignumber#1#2{%
  \ifx#2\end
    #1\let\next\@gobble
  \else
    #1\hspace{0pt plus 1pt}\let\next\@bignumber
  \fi
  \next#2}
\newcommand{\bignumber}[1]{\@bignumber#1\end}
\makeatother

\begin{document}

\date{\today}
\pagenumbering{arabic}

\maketitle

\begin{abstract}
\noindent 
For $q$ an odd prime power, we prove that there are always four  consecutive primitive elements in the finite field $\mathbb{F}_{q}$ when $q> 2401$.

 \end{abstract}
 

\section{Introduction}

Let $q= p^{n}$ for some $n\geq 1$ be the power of a prime $p$, and let $\mathbb{F}_{q}$  denote the finite field of size $q$, and $\mathbb{F}_{q}^{*}$ the set of non-zero elements in $\mathbb{F}_{q}$. 
An element $g \in \mathbb{F}_{q}$ is called a primitive element if it generates  $\mathbb{F}_{q}^{*}$. 

The problem of finding consecutive primitive elements in a finite field $\mathbb{F}_{q}$ is motivated by Brauer \cite{Brauer}, who, in 1928 examined long runs of consecutive quadratic non-residues. Brauer's result has been followed by the work of many authors (see, e.g., \cite{Hudson,Hummel}). Vegh \cite{Vegh} mentions a question raised by Brauer as to whether there are long runs of consecutive primitive roots modulo a prime $p$. Vegh proved this for pairs of consecutive primitive roots, seemingly unaware of a much more general result given earlier by Carlitz.

 Carlitz \cite{Carlitz} showed that, given any $n$, one may find a $q_{0}(n)$ such that $\mathbb{F}_{q}$ contains $n$ consecutive primitive elements for all $q> q_{0}(n)$. A natural question is: how does $q_{0}(n)$ grow with $n$? Cohen \cite{Cohen1,Cohen2,Cohen3} proved that $q_{0}(2) = 7$; Cohen, Oliveira e Silva and Trudgian \cite{COT2015} proved that $q_{0}(3) = 169$ --- see also \cite[Sect.\ 1]{COT2015} for a general bound on $q_{0}(n)$ and for more history on this problem. In both the $n=2$ and $n=3$ cases it is easy to provide a list of small values of $q$ for which $\mathbb{F}_{q}$ does not contain two or three primitive elements.

Turning to $n=4$, the best bound to date comes from Table 1 in \cite{COT2015}, namely that $q_{0}(4) \leq 3.29 \times 10^{32}$. In \cite[Sect.\ 6]{COT2015}, the authors conjecture that $q_{0}(4) = 7^{4} = 2401$. That is, for all $q>2401$ the finite field $\mathbb{F}_{q}$ should have four consecutive primitive elements. Numerical evidence in \cite{COT2015} shows that this is true for $2401<q<10^{8}$, whence `all' that we need to do is to check values of $q \in (10^{8}, 3.29 \times 10^{32})$. 

Clearly it is infeasible to check \textit{all} prime powers $q$ in this range. What is required is a method that allows us to get away without checking a large proportion of prime powers. Using a variant of the `prime divisor tree' (first announced in \cite{MTT} and further developed in \cite{Hunter} and \cite{JT}) we are able to resolve this conjecture.

\begin{theorem}
\label{conj:main1}
The finite field $\mathbb{F}_{q}$ has four consecutive primitive elements except when $q$
is divisible by $2$ or by $3$, or when q is one of the following: 
$5, 7, 11, 13, 17, 19, 23, 5^{2} , 29, 31, 41, 43,$ 
$61,67, 71, 73,79, 113, 11^{2} , 13^{2} , 181, 199,$
$337, 19^{2} , 397, 23^{2} , 571, 1093, 1381, 7^{4} = 2401$.
\end{theorem}

It was also conjectured in \cite[Sect.\ 6]{COT2015} that for all $q>15625$ the finite field $\mathbb{F}_{q}$ should have five consecutive primitive elements. We are not able to resolve this conjecture, but indicate, in Section \ref{S5}, some partial progress on it.
%
It does not seem feasible to resolve completely the problem of five consecutive primitive elements without either a new idea or a large increase in computational power.

This paper is organised as follows. In Section \ref{Outline} we outline the necessary background for the ensuing theory and computation. In Section \ref{S4} we list our algorithms and prove Theorem \ref{conj:main1}. In Section \ref{S5} we make some partial progress on the problem of five consecutive primitive roots.
Throughout this paper we let $\omega = \omega(q-1)$ denote the number of distinct prime factors of $q-1$, and let $\phi(n)$ denote Euler's totient function.

\section{Outline of the problem}\label{Outline}

We consider $\mathbb{F}_{q}$ separately, in the two cases when $q=p$
and when $q=p^{n}$, where $n \ge 2$. 
In the first case, to determine the primitive elements in $\mathbb{F}_{q}$ we use Pollard's factorisation \cite{JmP74}, that is, 
$q - 1 =  p^{a_{1}}_{1} \cdots p^{a_{r}}_{r}.$
Then $\alpha$ is a primitive element, if and only if 
$\alpha^{\frac{q-1}{p_{i}}} \neq 1,$ for all $1\leq i \leq r$. In the first case we use 
 Algorithm \ref{Algo2:Sieving} in section \ref{app:pdt}.

We outline below the procedure for the second case in which $q=p^{n}$ for $n \geq 2$. 
We use this in Algorithm \ref{Algo2:primitive} in Section \ref{Algae}.

We coded Algorithm \ref{Algo2:Sieving} using C/C++ GMP in  parallel multi-threading programming with OpenMP (Open Multi-processing) implementation for $q = p$
on desktop a PC with a quad core of 3.4 GHZ Intel i5 processors. 
An OpenMP is a library for parallel programming in the SMP (symmetric multi-processors, or shared-memory processors) model and all threads share memory and data.

The output of the results and the check for four consecutive primitive elements  
are given in Table~\ref{Table:tabe19}. 

On the same platform we coded  Algorithm \ref{Algo2:primitive} using \textit{Magma} for $q = p^{n}$, where $n \geq 2$, and the results are  given in Table~\ref{Table:tabe20}.

\subsection{Polynomial representation of primitive elements}

Recall that a monic irreducible polynomial whose roots are primitive elements in $\mathbb{F}_{q}$ is called a primitive polynomial.
It is well known that the field $\mathbb{F}_{q}$ can be constructed as $\mathbb{F}_{p}[x]/(f(x))$,
where $f(x)$ is an irreducible polynomial of degree $n$ over $\mathbb{F}_{p}$  and, in addition, if
$f(x)$ is primitive, then  $\mathbb{F}_{q}^{*}$ is generated multiplicatively by any root of $f(x)$. Note that $(f(x))$ is the maximal principal ideal generated by $f(x)$, i.e., it is an algebraic extension of   $\mathbb{F}_{p}$.
The field is exactly the set of all polynomials of degree $0$ to $n-1$ with the two field operations being addition and multiplication of polynomials modulo $f(x)$ and with modulo $p$ integer arithmetic on the polynomial coefficients.

We shall continually use the fact that
$a^{m} \in \mathbb{F}_{q}^{*}$ is a primitive element if and only if $(m, q - 1) = 1$, for some $m \in [1,\ldots,q-1]$.
We give two  examples of finite fields with  polynomial representations, the first with two, and the second with three consecutive primitive elements. It is straightforward to continue this to $N$ consecutive primitive elements.

\begin{enumerate}
	\item Consider the finite field $\mathbb{F}_{3^{2}} = \mathbb{F}_{3}[x]/(x^{2} + x +2)$, which is the set of polynomials 
	$$\{ 0, 1, 2, x, x+1, x+2, 2x, 2x+1, 2x+2 \},$$           with addition and multiplication of polynomials modulo $f(x) = x^{2} + x + 2$ and also modulo $3$. Note also that $f(x)$ is a primitive polynomial. 

There are $\phi(q-1)$ primitive elements,  that is  $\phi(8) = 4$ primitive elements, which are 
$\{a,a^{3},a^{5},a^{7}\} = \{2x,x+1,x,2x+2\}$, of which, clearly, $x$ and $x+1$ are consecutive.

\item Consider the finite field $\mathbb{F}_{7^{2}} = \mathbb{F}_{7}[x]/(x^2 + 6x + 3) = $
\begin{equation*}
\begin{split}
\{&0, 1, 2, 3, 4, 5, 6, x,x + 1,x + 2,x + 3, x + 4, x + 5,x + 6, 2x, 2x + 1, 2x + 2, 2x + 3,\\ 
& 2x + 4, 2x + 5, 2x + 6, 3x, 3x + 1, 3x + 2, 3x + 3, 3x + 4, 3x + 5, 3x + 6, 4x, 4x +1\\
& 4x + 2, 4x + 3, 4x + 4, 4x + 5,4x + 6, 5x, 5x + 1, 5x + 2, 5x + 3, 5x + 4, 5x + 5\\
& 5x + 6, 6x, 6x + 1, 6x + 2, 6x + 3, 6x + 4, 6x + 5, 6x + 6\},
\end{split}
\end{equation*}
where $f(x) = x^2 + 6x + 3$ is a primitive polynomial.
There are $\phi(q-1)$ primitive elements in $\mathbb{F}_{7^{2}}$, that is  $\phi(48) = 16$ primitive elements, and these are: 
\begin{equation*}
\begin{split}
\{&x, x + 1, x + 5, x + 6, 2x, 2x + 5, 3x + 1, 3x + 3,\\
 &4x + 4, 4x + 6,
5x, 5x + 2, 6x, 6x + 1, 6x + 2, 6x + 6 \}.
\end{split}
\end{equation*}
Three consecutive primitive elements are $6x,  6x + 1, 6x + 2 $; note also that $x+5, x+6, x$ are also three consecutive primitive elements. Throughout the remainder of this paper we are never concerned with multiple sets of $N$ consecutive roots: we merely wish to check (quickly) that there is at least one such set.

\end{enumerate}

\subsection{Sieving preliminaries}\label{prim}
We use the following criterion to prove the existence of $n$ consecutive primitive elements in $\mathbb{F}_{q}^{*}$, for $n=4,5$. 

\begin{lemma}[Theorem 5 \cite{COT2015}]
\label{lem:main1}
Suppose $3 \leq n \leq p$ and $e$ is a divisor of $q-1$. If 
$Rad(e) = Rad(q - 1)$, then
set $s = 0$ and $\delta = 1$. Otherwise, let $p_{1}, \cdots, p_{s}$ , $s \geq 1$, 
be the primes dividing $q - 1$ but not $e$
and set $\delta = 1 - n \sum_{i=1}^{s} p_{i}^{-1}$. Assume that $\delta > 0$. If also 
\begin{align}
\label{sievingInq}
q > \left((n-1)\left(\frac{ns-1}{\delta} + 2\right)(2^{n(\omega-s)}) \right)^{2},
\end{align}
then there exist $n$ consecutive primitive elements in $\F_{q}$.
\end{lemma}

We have another, better, criterion when $q \equiv 3\pmod{4}$ as follows.
\begin{lemma}[Theorem 6 \cite{COT2015}]
\label{lem:main2}
Suppose that $q\equiv 3 \pmod{4}$ and that $3 \leq n \leq p$ and $e$ is an even divisor of $q -1$. If 
$Rad(e) = Rad(q - 1)$, then
set $s = 0$ and $\delta = 1$. Otherwise, let $p_{1}, \cdots, p_{s}$ , $s \geq 1$, 
be the primes dividing $q - 1$ but not $e$
and set $\delta = 1 - n \sum_{i=1}^{s} p_{i}^{-1}$. Assume that $\delta > 0$. If also 
\begin{align}
\label{sievingInq2}
q > \left(\frac{(n-1)}{2}\left(\frac{ns-1}{\delta} + 2\right)(2^{n(\omega-s)}) \right)^{2},
\end{align}
then there exist $n$ consecutive primitive elements in $\F_{q}$.
\end{lemma}
When $q\equiv 3\pmod{4}$ we use (\ref{sievingInq2}), which gives an improvement of a factor of $4$ over the criterion in (\ref{sievingInq}). When $q\equiv 1\pmod{4}$, while we are forced to use the inferior bound in (\ref{sievingInq}), we are still able to obtain a small improvement. We know in this case that $2^{2}|q-1$. This helps in constructing possible counterexamples that require further checking. It is the use of this `divide and conquer' approach along with the implementation of the prime divisor tree that enables us to prove Theorem \ref{conj:main1}.

We conclude this section by mentioning two relevant results from Table 1 in \cite{COT2015}. For the problem of four consecutive primitive elements we need only consider those $q\leq 3.29\times 10^{32}$ such that $\omega(q-1) \leq 23$. For five consecutive primitive elements we need only examine those $q\leq 4.22\times 10^{61}$ with $\omega(q-1) \leq 37$.

\section{Four consecutive primitive elements}\label{S4}

We first present Algorithm \ref{Algo1}, which determines the choice of $s$ and $\delta$ to minimise the right hand sides of~\eqref{sievingInq} and~of~\eqref{sievingInq2}. 
First, set integral intervals for $s  \in [a,b]$ and $\omega \in [b,c]$ such that $a < b < c$.

\begin{algorithm}[H]
	\DontPrintSemicolon
	\KwIn{$a,b,c,n:$ { $s \in[a,b]$ and $w \in [b,c]$, where $n$ is number of consecutive primitive roots.}} 
	\KwResult{$M$} 
	\Fn{SievingAlgorithm(a,b,c)}{
		Let $w = \omega(q-1)$\;
		$M = \{\}$\;   
		\For{$s \in [a,b]$}{ 
			\For{$w \in [b,c]$}{
				$L=[2,3,5,7, \ldots , w=p_{\omega(q-1)}]$ {\tiny\tcp*[h]{\textbf{list of distinct primes}}}.\;
				$w_{e} = w-s$\; 
				Assert   {$ w_{e} \geq 0$}\;
				$L_0=[2,3, \ldots , w_e=q_{\omega(q-1)}]$ {\tiny\tcp*[h]{\textbf{list of  $e \mid q-1$}}}\;  
				$\widetilde{L} = (set(L) - set(L_0))$ {\tiny\tcp*[h]{\textbf{remove $e \mid p-1$ from  the list $L$}}}\;
				$d = \sum_{p \in \widetilde{L}} 1/p$\;
				$\delta = 1- n*d$\;
				\If{ $\delta > 0$} {
				   \eIf { $q \equiv 1 \; \pmod{4}$} {
					  Evaluate the sieving inequality~\eqref{sievingInq}
					   \begin{align*} 
					    R &= \prod_{p \in L}p\\
					    S &= \left((n-1)\left(\frac{ns-1}{\delta} + 2\right)(2^{n(w-s)})\right)^{2}
					    \end{align*}   
					  }{	    
				       Evaluate the sieving inequality~\eqref{sievingInq2}
					    \begin{align*} 
					    R &= \prod_{p \in L}p\\
					    S &= \left(\frac{(n-1)}{2}\left(\frac{ns-1}{\delta} + 2\right)(2^{n(w-s)})\right)^{2}
					   \end{align*}  
					}  
				     \If{$R > S$}{
						 append $M[s] = w, \delta$\;   
					 }	
				} 
			} 
		} 
		\Return{$M$}\;
	} 
	\caption{\small \bf{Determining the range of $w$ for choices of $s$}}
	\label{Algo1}
\end{algorithm}

For each $s \in [a,b]$ we find values of  $\omega \in [b,c]$ that satisfy the sieving  inequalities~\eqref{sievingInq} and~\eqref{sievingInq2}, provided that $\delta > 0$, then for each value of $s$ append ranges of $\omega$ into a list.

For each value of $3 \le \omega(q-1) \le 23$, 
we applied Lemmas \ref{lem:main1} and \ref{lem:main2} directly. Algorithm~\ref{Algo1} was then used to generate the number of possible exceptions 
to Theorem \ref{conj:main1}, which are in the fifth column of Table \ref{Table_n4}. 
There are clearly too many possible exceptions to check. To resolve this we now define the  prime divisor tree.

\begin{table}[H]
	\centering
	\begin{tabular}{c c c c c }
		\hline\hline	
		\multicolumn{5}{c}{{\bf Values of $\omega(q-1)$ and total exceptions for four consecutive primitive elements}}\\ \hline
		$\omega(p-1)$ & $s$  & $\delta$   &Intervals & Possible exceptions\\ \hline
		$23$          &$12$ & $0.1376994749839410$ &$(2.670  \times 10^{32},   5.580  \times  10^{33})$   & $2.656 \times 10^{33}$  \\
		$22$          &$12$ & $0.0568599880037627$ &$(3.217  \times 10^{30},   3.790  \times  10^{31})$   & $1.734 \times 10^{31}$ \\
		$21$          &$11$ & $0.1074928993961680$ &$(4.072  \times 10^{28},   4.396  \times  10^{29})$   & $1.994 \times 10^{29}$ \\ 
		$20$          &$11$ & $0.0243563854613543$ &$(5.579  \times 10^{26},   3.319  \times  10^{28})$   & $1.631 \times 10^{28}$ \\
		$19$          &$10$ & $0.0806944136303683$ &$(7.858  \times 10^{24},   2.502  \times  10^{27})$   & $1.247 \times 10^{27}$ \\
		$18$          &$9$  & $0.1403959061676820$ &$(1.172  \times 10^{23},   6.709  \times  10^{26})$   & $3.354 \times 10^{26}$ \\
		$17$          &$9$  & $0.0320566331812242$ &$(1.922  \times 10^{21},   4.965  \times  10^{25})$   & $2.482 \times 10^{25}$ \\ 
		$16$          &$8$  & $0.0998532433507158$ &$(3.258  \times 10^{19},   4.052  \times  10^{24})$   & $2.026 \times 10^{24}$ \\
		$15$          &$7$  & $0.1753249414639230$ &$(6.148  \times 10^{17},   1.010  \times  10^{24})$   & $5.050 \times 10^{23}$ \\ 
		$14$          &$7$  & $0.0499050086531730$ &$(1.308  \times 10^{16},   4.780  \times  10^{22})$   & $2.390 \times 10^{22}$ \\ 
		$13$          &$6$  & $0.1429282644671260$ &$(3.042  \times 10^{14},   4.303  \times  10^{21})$   & $2.151 \times 10^{21}$ \\ 
		$12$          &$5$  & $0.2404892400768830$ &$(7.420  \times 10^{12},   1.063  \times  10^{21})$   & $5.319 \times 10^{20}$ \\ 
		$11$          &$5$  & $0.1133032305379320$ &$(2.005  \times 10^{11},   1.823  \times  10^{19})$   & $9.118 \times 10^{18}$ \\ 
		$10$          &$4$  & $0.2423354886024480$ &$(6.469  \times 10^{9},    2.585  \times  10^{18})$   & $1.292 \times 10^{18}$ \\ 
		$9$           &$4$  & $0.0725742153928989$ &$(2.230  \times 10^{8},    1.077  \times  10^{17})$   & $5.386 \times 10^{16}$ \\ 
		$8$           &$3$  & $0.2464872588711600$ &$(9.699  \times 10^{6},    5.378  \times  10^{15})$   & $2.689 \times 10^{15}$\\ 
		$7$           &$3$  & $0.0933772110242699$ & $(5.105  \times 10^{5},   1.386  \times  10^{14})$   & $6.098 \times 10^{13}$ \\ 
		$6$           &$2$  & $0.3286713286713290$ & $(3.003  \times 10^{4},   5.245  \times  10^{12})$   & $2.622 \times 10^{12}$ \\ 
		$5$           &$1$  & $0.6363636363636360$ & $(2.311  \times 10^{3},   4.551  \times  10^{11})$   & $2.178 \times 10^{11}$ \\ 
		$4$           &$1$  & $0.4285714285714290$ &$(2.110  \times 10^{2},    3.057  \times  10^{9})$    & $1.528 \times 10^{9}$ \\ 
		$3$           &$1$  & $0.2000000000000000$ &$(3.100  \times 10^{1},    4.261  \times  10^{7}) $   & $1.887 \times 10^{7}$\\ 	
		\hline\hline
	\end{tabular} 
	\caption{Choices of $s$  and $\delta$ for values of $\omega(q-1)$ for four consecutive primitive elements.}
	\label{Table_n4}
\end{table}

\subsection{The prime divisor tree }
\label{app:pdt}
The point of the algorithm is to split the problem into many sub-cases, according as $p|q-1$ or not, where $p$ is a `small' prime. Since the size of $\delta$ in Lemmas \ref{lem:main1} and \ref{lem:main2} depends precisely on small prime factors of $q-1$, this approach allows for more specific information to be wrought from the sieve.

\begin{algorithm}[H]
	\DontPrintSemicolon
	\KwData{$L=[2,3,5,7, \ldots ,q_{w}]$ list of distinct primes, let $w=\omega(q-1)$, where $n$ is number of consecutive integers.}
	\KwIn{An interval  $I = (lower, upper)$ see Table~\ref{Tab:Tablecomp}.}
	\KwResult{$D_{s} = \prod_{p \in M[s]}p$, where $p \mid q-1$, with respect to $s$.}
	
	\Fn{PrimeDivisorTree(a,w,n)}{
		$M = \{ \}$\;  
		\For{$i=0; i< size(L); i= i+1$}{
			let $t = L[i]$ and assume that $t \nmid p-1$\;
			$\widetilde{L} = (L - set(t))$  {\tiny\tcp*[h]{\textbf{remove $t$ from  the list $L$}}}\;
			$x = q_{w+1}$ {\tiny\tcp*[h]{\textbf{get new prime the $(w+1)$th prime to replace $t$.}}}\;
			append $x$ to $\widetilde{L}$\;
			\For{$s \in [a,w]$}{
				$w_e = w-s$\;
				\If{ $w_e < 0$} {
					continue\;
				}
				$L_{s} = \widetilde{L}[w_{e}:]$ {\tiny\tcp*[h]{\textbf{remove $w_{e}$ elements from $L$ starting from the 
							begining index.}}}\; 
				$d = \sum_{p \in L_{s}} 1/p$\;
				$\delta = 1- n*d$\;
				\If{ $\delta > 0$} {
				
				\eIf { $q \equiv 1 \; \pmod{4}$} {
					  Evaluate the sieving inequality~\eqref{sievingInq}
					   \begin{align*} 
					    R &= \prod_{p \in L}p\\
					    S &= \left((n-1)\left(\frac{ns-1}{\delta} + 2\right)(2^{n(w-s)})\right)^{2}
					    \end{align*}   
					  }{	    
				       Evaluate the sieving inequality~\eqref{sievingInq2}
					    \begin{align*} 
					    R &= \prod_{p \in L}p\\
					    S &= \left(\frac{(n-1)}{2}\left(\frac{ns-1}{\delta} + 2\right)(2^{n(w-s)})\right)^{2}
					   \end{align*}  
					}  
				     \If{$R > S$}{
						 append $M[s] = w, \delta$\;   
					 }	
				} 
			} 
			$D_{s} = \prod_{p \in M}p$ {\tiny\tcp*[h]{\textbf{product of  $p\in M$ where  $p \mid q-1$}}}.\;
			\Return{$D_{s}$}\;
		} 
	} 
	\caption{\sc{Prime divisor} tree}
	\label{Algo:level1}
\end{algorithm}
\vspace{0.05cm}
We note that we can keep splitting into further sub-cases if required. 
When we have $k$ such cases (that is, conditions on the first $k$ primes dividing, or not dividing, $q-1$) we say that we have gone down the prime divisor tree to level~$k$.

Using Algorithm~\ref{Algo:level1}, we see that the choice of $s=12$ and $\delta =0.137699474983941$ eliminates $\omega = 23$ immediately. Similarly choosing $s=12$ and $\delta =0.0568599880037627$ eliminates $\omega = 22$.
This results in the values summarised in Table~\ref{Tab:Tablecomp}.

\vspace{0.3cm}
\begin{table}[H]
	\centering
	\begin{tabular}{c c c c c c}
		\hline\hline
		$\omega(q-1)$ &s  &$\delta$ & \; \; \; \; \; Interval &Possible             &Prime divisor\\
		              &    &        &                         &exceptions           &tree level\\ [0.5ex]\hline	              
		  $21$          &$11$ & $0.1074928993961680$ &$(4.072  \times 10^{28},   4.396  \times  10^{29})$   & $29$                   &$1$\\ 
		  $20$          &$11$ & $0.0243563854613543$ &$(5.579  \times 10^{26},   3.319  \times  10^{28})$   & $175$                  &$1$ \\
		  $19$          &$10$ & $0.0806944136303683$ &$(7.858  \times 10^{24},   2.502  \times  10^{27})$   & $952$                  &$1$\\
		  $18$          &$9$  & $0.1403959061676820$ &$(1.172  \times 10^{23},   6.709  \times  10^{26})$   & $85796$                &$2$\\
		  $17$          &$9$  & $0.0320566331812242$ &$(1.922  \times 10^{21},   4.965  \times  10^{25})$   & $387383$               &$2$\\ 
		  $16$          &$8$  & $0.0998532433507158$ &$(3.258  \times 10^{19},   4.052  \times  10^{24})$   & $1.865 \times 10^{6}$  &$2$\\
		  $15$          &$7$  & $0.1753249414639230$ &$(6.148  \times 10^{17},   1.010  \times  10^{24})$   & $1.724 \times 10^{8}$  &$3$\\ 
		  $14$          &$7$  & $0.0499050086531730$ &$(1.308  \times 10^{16},   4.780  \times  10^{22})$   & $3.837 \times 10^{8}$  &$3$\\ 
		  $13$          &$6$  & $0.1429282644671260$ &$(3.042  \times 10^{14},   4.303  \times  10^{21})$   & $1.485 \times 10^{9}$  &$3$ \\ 
		  $12$          &$5$  & $0.2404892400768830$ &$(7.420  \times 10^{12},   1.063  \times  10^{21})$   & $1.655 \times 10^{11}$ &$4$\\ 
		  $11$          &$5$  & $0.1133032305379320$ &$(2.005  \times 10^{11},   1.823  \times  10^{19})$   & $1.050 \times 10^{11}$ &$4$\\ 
		  $10$          &$4$  & $0.2423354886024480$ &$(6.469  \times 10^{9},    2.585  \times  10^{18})$   & $2.308 \times 10^{11}$ &$4$\\ 
		  $9$           &$4$  & $0.0725742153928989$ &$(2.230  \times 10^{8},    1.077  \times  10^{17})$   & $2.788 \times 10^{11}$ &$4$\\ 
		  $8$           &$3$  & $0.2464872588711600$ &$(9.699  \times 10^{6},    5.378  \times  10^{15})$   & $3.202 \times 10^{11}$ &$4$\\
		  $7$           &$3$  & $0.0933772110242699$ &$(5.105  \times 10^{5},    1.386  \times  10^{14})$   & $2.039 \times 10^{12}$ &$5$\\ 
		  $6$           &$2$  & $0.3286713286713290$ &$(3.003  \times 10^{4},    5.245  \times  10^{12})$   & $2.622 \times 10^{12}$ &$0$\\ 
		  $5$           &$1$  & $0.6363636363636360$ &$(2.311  \times 10^{3},    4.551  \times  10^{11})$   & $2.178 \times 10^{11}$ &$0$ \\ 
		  $4$           &$1$  & $0.4285714285714290$ &$(2.110  \times 10^{2},    3.057  \times  10^{9})$    & $1.528 \times 10^{9}$  &$0$\\ 
		  $3$           &$1$  & $0.2000000000000000$ &$(3.100  \times 10^{1},    4.261  \times  10^{7})$    & $1.887 \times 10^{7}$  &$0$\\ 	
		  \hline\hline            
	\end{tabular}
	\caption{Intervals containing exceptions to
		Theorem~\ref{conj:main1} for a given value of $\omega(q-1)$ with prime divisor tree levels.}
	\label{Tab:Tablecomp}
\end{table}
We now present Algorithm \ref{Algo2:Sieving} that allows us to check, from our list of possible exceptions in Table \ref{Tab:Tablecomp}, whether or not $q$ has four consecutive primitive roots.

\vspace{0.3cm}
\begin{algorithm}[H]
	\DontPrintSemicolon
	\KwData{Interval $I =(lower, upper)$ in Table~\ref{Tab:Tablecomp} }
	\KwIn{ $D = \prod p$, where $ p\mid q-1$ from Algorithm~\ref{Algo:level1} }
	\KwResult{Return initial list of primes for interval $I$}
	\Fn{Sieving algorithm}{
		Find initial positive integer $m$ such that $D \mid m$ and $m := min(D\Z \cap I)$\;
		$S := [ \; ]$\;
		set $w \in \{23,22,21, 20, 19, 18, 17, 16, 15, 14, 13, 12, 11, 10, 9, 8, 7, 6, 5, 4, 3\}$\;
		\For{$n = m$; $\;$ $n \le upper$; $\;$ $n = n+D_{w}$}{
			Assert   {$n \% D == 0$}\;
			$p    = n + 1$\;  
			\If{$Isprime(q)$}{ 
				\If{$\omega(q-1)== w$}{     
					append $q$ to $S$ {\tiny\tcp*[h]{\textbf{save the initial list of primes.}}}\;
					Find $4$ consecutive primitive elements $\pmod{q}$\;
				}          
			}
		} 
		\Return{$S$}\; 
	} 
	\caption{Sieving for initial list of primes}
	\label{Algo2:Sieving}
\end{algorithm}

\vspace{0.3cm}

We include here an example of Algorithm~\ref{Algo2:Sieving}. Consider one of the possible 29 exceptions for $w=21$. We have that $D$ must divide $q-1$, where
$$D= \prod_{p \in M}p =13576560199749674716873774490.$$ The choice of $s=11$ gives 
\begin{equation}\label{light}
q-1 = D*m \in (4.072 \times 10^{28}, 4.396 \times 10^{29}),
\end{equation} for some integer $m>0$.
There is only one prime $q$ satisfying (\ref{light}). This corresponds to $m = 18$, namely, 
$$ q = D*m + 1 = 244378083595494144903727940821,$$ and its four consecutive primitive roots are $993, 994, 995, 996$. 

We repeat the procedure for $3 \leq \omega \leq 20$ to verify Theorem~\ref{conj:main1}, in the case of $q$ prime. The output is summarised in Table~\ref{Table:tabe19}.

\begin{table}[H]\small
	\centering
	\begin{tabular}
		{l l l l l l l}
		\hline\hline
		\multicolumn{7}{c}{\bf Prime dividing $q-1$ using prime divisor tree with $7 \le \omega \le 21$}\\ \hline
		$\omega$ &$s$ & $ q\nmid (p-1)$	 &  $M=L\backslash\{p\}$   & $D = \prod_{q \in M}q$ & Number of  &Number of\\
	    	&    &                   &                         &                        & exceptions     &primes\\ \hline      
		$21$     &$11$  &$3$           &   $2, 5, 7, \cdots, 73$    &$1.3576 \times 10^{28}$ &$29$                    &$1$ \\ 
		$20$     &$11$  &$3$           &   $2, 5, 7, \cdots, 71$    &$1.8598 \times 10^{26}$ &$175$                   &$14$ \\
		$19$     &$11$  &$3,5$         &   $2, 7, 11, \cdots, 67$   &$5.2388 \times 10^{23}$ &$952$                   &$44$ \\ 
		$18$     &$10$  &$3,5$         &   $2, 7, 11, \cdots, 61$   &$7.8192 \times 10^{21}$ &$85796$                 &$2794$ \\ 
		$17$     &$10$  &$3,5$         &   $2, 7, 11, \cdots, 59$   &$1.2818 \times 10^{20}$ &$387383$                &$13255$ \\ 
		$16$     &$9$   &$3,5,7$       &   $2, 11, 13, \cdots, 53$  &$3.1037 \times 10^{17}$ &$1.865 \times 10^{6}$   &$66420$ \\ 
		$15$     &$9$   &$3,5,7$       &   $2, 11, 13, \cdots, 47$  &$5.8560 \times 10^{15}$ &$1.724 \times 10^{8}$   &$3566570$ \\ 
		$14$     &$9$   &$3,5,7$       &   $2, 11, 13, \cdots, 41$  &$2.8976 \times 10^{12}$ &$3.837 \times 10^{8}$   &$8464555$ \\ 
		$13$     &$8$   &$3,5,7,11$    &   $2, 13, 17, \cdots, 41$  &$2.6342 \times 10^{11}$ &$1.485 \times 10^{9}$   &$34451458$ \\       
		$12$     &$8$  &$3,5,7,11$     &   $2, 13, 17, \cdots, 37$  &$6424881502$            &$1.655 \times 10^{11}$  &$2063206920$ \\ 
		$11$     &$7$  &$3,5,7,11$     &   $2, 13, 17, \cdots, 31$  &$173645446$             &$1.050 \times 10^{11}$  &$1405800079$ \\ 
	    $10$     &$7$  &$3,5,7,11$     &   $2, 13, 17, \cdots, 29$  &$11202932$              &$2.308 \times 10^{11}$  &$3291590809$ \\ 
	    $9$      &$6$  &$3,5,7,11,13$  &   $2, 17, 19, 23$          &$386308$                &$2.788 \times 10^{11}$  &$4293202707$ \\ 
	    $8$      &$6$  &$3,5,7,11,13$  &   $2, 17, 19$              &$16796$                 &$3.202 \times 10^{11}$  &$5239259761$ \\ 
	    $7$      &$6$  &$3,5,7,11,13$  &   $2, 17$                  &$68$                    &$2.039 \times 10^{12}$  &$15550071093$ \\ 
		$6$     &$2$   &$  $           &                            &                        &$1.311 \times 10^{12}$  &$5480805894$ \\ 
		$5$     &$1$   &$  $           &                            &                        &$2.178 \times 10^{11}$  &$2250816606$ \\ 
		$4$     &$1$   &$  $           &                            &                        &$1.528 \times 10^{9}$   &$27237430$\\
		$3$     &$0$   &$  $           &                            &                        &$1.887 \times 10^{7}$   &$382271$ \\
		\hline\hline
	\end{tabular}
	\caption{Number of primes with $3 \le \omega \le 21$ and number of exceptions.}
	\label{Table:tabe19}
\end{table}


\subsection{Four consecutive primitive elements in $\mathbb{F}_{p^{n}}$, with $n \ge 2$}\label{Algae}
In this section we adopt the following procedure when $q=p^{n}$ for $n\geq 2$. 

\begin{enumerate}
	\item Find a primitive polynomial $f(x)$ with degree $n$. 
	\item Construct $\mathbb{F}_{q} =  \mathbb{F}_{p}{[x]}/(f(x))$.
		\item If $x \in  \mathbb{F}_{q}^{*}$, then $x^{m}$ is primitive if and only if $(m, q-1) = 1$, for some $m \in [1, \ldots, q-1]$.
	\item Find consecutive primitive elements among $x^{m}\pmod{f(x)}$.
\end{enumerate}
This is summarised in Algorithm \ref{Algo2:primitive} below.

\vspace{0.3cm}

\begin{algorithm}[H]
	\DontPrintSemicolon
	\KwData{Given a prime power $q= p^{n}$, $n \ge 2$}
	\KwIn{ Finite field Prime $\mathbb{F}_{q}^{*}$ }
	\KwResult{Return primitive element modulo $f(x)$, where $f(x)$ is a primitive polynomial}
	\Fn{PrimitiveElement()}{
		generate randomly $\alpha \in \mathbb{F}_{q}^{*}$\;
		\For{each $i \in [1..q-1]$ }{
			\If{$GCD(i, q-1)==1$ }{
				then $\alpha^{i}$ {\textbf is a primitive element}\;
			}    
			\Return{$\alpha^{i} \pmod{f(x)}$}\;   
		}            
	} 
	\caption{Compute primitive element modulo $f(x)$}
	\label{Algo2:primitive}
\end{algorithm}

\vspace{0.3cm}

In Table \ref{Table:tabe20} we give examples of four consecutive primitive elements in $\mathbb{F}_{q}[x]$ for $q =p^{2}$.
	\begin{table}[H]\small
		\centering
		\begin{tabular}
			{l c c c c}
			\hline\hline
			\multicolumn{5}{c}{{\bf Four consecutive primitive elements in  $\mathbb{F}_{q}$ where $q =p^{2}$}}\\ \hline  
			$\omega$     &$p$            &$q = p^{2}$     &Four consecutive primitive elements & Primitive polynomials \\
			4      &$29$            & $841$          & $11x$, $11x + 1$, $11x + 2$, $11x + 3$, & $f(x)= x^2 + 24x + 2$\\
			3      &$31$            & $961$          & $10x +4$, $10x + 5$, $10x + 6$, $10x + 7$, & $f(x)= x^2 + 29x + 3$\\
			3      &$37$            & $1369$         & $2x + 13$, $2x + 14$, $2x + 15$, $2x + 16$, & $f(x)= x^2 + 33x + 2$\\
			4      &$41$            & $1681$         & $15x + 16$, $15x + 17$, $15x + 18$, $15x + 19$, & $f(x)= x^2 + 38x + 6$\\
			4      &$43$            & $1849$         & $5x + 39$, $5x + 40$, $5x + 41$, $5x + 42$, & $f(x)= x^2 + 42x + 3$\\			
			3      &$47$            & $2209$         & $7x + 14$, $7x + 15$, $7x + 16$ $7x + 17$, & $f(x)=x^2 + 45x + 5$\\
			3      &$53$            & $2809$         &  $2x + 3$, $2x + 4$, $2x + 5$, $2x + 6$, & $f(x)=x^2 + 49x + 2$\\ 
			4      &$59$            & $3481$         & $11x + 46$, $11x + 47$, $11x + 48$, $11x + 49$, &$f(x)=x^2 + 58x + 2$\\
			4      &$61$            & $3721$         & $5x + 9$, $5x + 10$, $5x + 11$, $5x + 12$, &$f(x)=x^2 + 60x + 2$\\
			4      &$67$            & $4489$         & $2x + 13$, $2x + 14$, $2x + 15$, $2x + 16$, &$f(x)=x^2 + 63x + 2$\\
			4      &$71$            & $5041$         & $22x + 12$, $22x + 13$, $22x + 14$, $22x + 15$, &$f(x)=x^2 + 69x + 7$\\
			3      &$73$            & $5229$         & $x + 8$, $x + 9$, $x + 10$, $x + 11$, &$f(x)=x^2 + 70x + 5$\\
			4      &$79$            & $6241$         & $15x + 62$, $15x + 63$, $15x + 64$, $15x + 65$, &$f(x)=x^2 + 78x + 3$\\
			4      &$83$            & $6889$         & $8x + 18$, $8x + 19$, $8x + 20$, $8x +21$, & $f(x)=x^2 + 82x + 2$\\
			4      &$89$            & $7921$         & $24x + 29$, $24x + 30$, $24x + 31$, $24x + 32$, &$f(x)=x^2 + 82x + 3$\\	 
			3      &$97$            & $9409$         & $6x + 3$, $6x + 4$, $6x + 5$, $6x + 6$, &$f(x)=x^2 + 96x + 5$\\
			4      &$101$           & $10201$        & $5x + 39$, $5x + 40$, $5x + 41$, $5x + 42$, &$f(x)=x^2 + 97x + 2$\\
			4      &$103$           & $10609$        & $2x + 43$, $2x + 44$, $2x + 45$, $2x + 46$, &$f(x)=x^2 + 102x + 5$\\	
		
			\vdots  &\vdots          & \vdots            & \vdots                              & \vdots \\
			
			5       &$131$           & $17161$           & $3x + 21$, $3x + 22$, $3x + 23$, $3x + 24$, &$f(x)=x^2 + 127x + 2$\\
			4       &$137$           & $18769$           & $x + 38$, $x + 39$, $x + 40$, $x + 41$, &$f(x)=x^2 + 131x + 3$\\
			5       &$139$           & $19321$           & $x + 30$, $x + 31$, $x + 32$, $x + 33$, &$f(x)=x^2 + 138x + 2$\\
			4       &$149$           & $22201$           & $2x + 26$, $2x + 27$, $2x + 28$, $2x + 29$, &$f(x)=x^2 + 145x + 2$\\
			
			\vdots  &\vdots         &  \vdots             &\vdots                              & \vdots  \\
			
			\hline\hline
		\end{tabular}
		\caption{ Four consecutive primitive elements in $\mathbb{F}_{q}[x]$ modulo $f(x)$.}
		\label{Table:tabe20}
	\end{table}

We continue by repeating the procedure of Algorithm \ref{Algo2:primitive} to verify the existence of four consecutive primitive elements of $\mathbb{F}_{q}[x]$ modulo primitive polynomials $f(x)$ for $q = p^{n}$, where $n \geq 2$. For example, the running time for $\omega = 3, 4, 5$ and $q = p^{2}$  using {\em Magma}~\cite{BCFS20} in  Table ~\ref{Table:tabe20} took $0.35$ seconds to generate the results.
For $\omega = 7$, in Table ~\ref{Table:tabe19}, it took approximately three weeks to check the $2.039 \times 10^{12}$ possible exceptions. For the remaining values of $\omega$ the running times are much less. This completes the computational part of proof of Theorem~\ref{conj:main1}.

\section{Conclusion: five consecutive primitive elements}\label{S5}
Since we have resolved the case of four consecutive primitive elements, it is natural to ask whether we can tackle five. In this section we show that this problem is out of reach for $5 \le \omega \le 25$ with current methods. 

We have included the results obtained by repeating the procedure of the preceding sections for five consecutive primitive elements with $26 \le \omega \le 37$. Using Lemmas \ref{lem:main1} and  \ref{lem:main2} the values $\omega = 36,37$ are eliminated immediately. 

Next we apply Algorithm~\ref{Algo1} on on $26 \le \omega(q-1) \le 35$, yielding Table~\ref{Tab:TableFiveCons}.

\begin{table}[H]
	\centering
	\begin{tabular}{c c c c c}
		\hline\hline
		\multicolumn{5}{c}{{\bf Values of $\omega$ and $s$ for four consecutive primitive elements  }}\\ \hline
		$\omega$ & $s$  & $\delta$ &Interval &Exceptions \\ \hline
        $35$          &$18$  & $0.05477236812843190$    &$(1.492  \times 10^{57},  1.584   \times  10^{58})$ & $1.43480 \times 10^{58}$ \\ 
		$34$          &$18$  & $0.00358365239643321$    &$(1.001  \times 10^{55},  3.606   \times  10^{57})$ & $3.59599 \times 10^{57}$\\ 
		$33$          &$17$  & $0.03955487541801610$    &$(7.205  \times 10^{52},  2.641   \times  10^{55})$ & $2.63379 \times 10^{55}$\\ 
		$32$          &$16$  & $0.07605122578297950$    &$(5.259  \times 10^{50},  6.332   \times  10^{54})$ & $6.33147 \times 10^{54}$\\ 
		$31$          &$16$  & $0.01987954207276780$    &$(4.014  \times 10^{48},  9.025   \times  10^{52})$ & $9.02459 \times 10^{52}$\\ 
		$30$          &$15$  & $0.05924962081292530$    &$(3.161  \times 10^{46},  8.934   \times  10^{51})$ & $8.93396 \times 10^{51}$\\ 
		$29$          &$14$  & $0.10349740842354500$    &$(2.797  \times 10^{44},  2.553   \times  10^{51})$ & $2.55299 \times 10^{51}$\\ 
		$28$          &$14$  & $0.04298598933316810$    &$(2.566  \times 10^{42},  1.440   \times  10^{49})$ & $1.43999 \times 10^{49}$\\ 
		$27$          &$13$  & $0.08971496129578480$    &$(2.398  \times 10^{40},  2.853   \times  10^{48})$ & $2.85299 \times 10^{48}$\\ 
		$26$          &$13$  & $0.02197958084873130$    &$(2.329  \times 10^{38},  4.623   \times  10^{46})$ & $4.62299 \times 10^{46}$\\ 
		\hline\hline
	\end{tabular}
	\caption{Intervals containing exceptions to five consecutive primitive elements.}
	\label{Tab:TableFiveCons}
\end{table}

 As before, the number of possible exceptions is reduced substantially with the prime divisor tree. This results in a small list of primes for further checking. We establish that each of these primes has five consecutive primitive elements.
 The output is summarised in Table~\ref{Table:tab16_4}.

\begin{table}[H]\small
	\centering
	\begin{tabular}
		{l l l l l l l}
		\hline\hline
		\multicolumn{7}{c}{\bf Prime dividing $q-1$ using prime divisor tree with $26 \le \omega \le 34$}\\ \hline
		$\omega$ &$s$ & $ q\nmid (p-1)$	 &  $M=L\backslash\{p\}$   & $D = \prod_{q \in M}q$ & Number of  &Number of\\
		         &       &               &                         &                            & exceptions &primes\\ \hline  
		$35$     &$18$   & $3$           &   $2,5,7, \cdots, 149$ &$4.973411 \times  10^{56}$    & $28$                      &$0$ \\      
		$34$     &$18$   & $3$           &   $2,5,7, \cdots, 139$ &$3.338215 \times  10^{54}$    & $1077$                    &$27$ \\
		$33$     &$18$   & $3$           &   $2,5,7, \cdots, 137$ &$2.401593 \times  10^{52}$    & $1096$                    &$37$ \\ 
		$32$     &$17$   & $3$           &   $2,5,7, \cdots, 131$ &$1.752988 \times  10^{50}$    & $36120$                   &$957$ \\ 	
	    $31$     &$16$   &$3,5$          &   $2,7,11, \cdots, 127$ &$ 2.676317 \times  10^{47}$  & $337195$                  &$6118$ \\   
	    $30$     &$16$   &$3,5$          &   $2,7,11, \cdots, 113$ &$2.107336 \times  10^{45}$   & $4.2394 \times 10^{6}$    &$83823$ \\ 
	    $29$     &$16$   &$3,5$          &   $2,7,11, \cdots, 109$ &$1.864899 \times  10^{43}$   & $1.36882 \times 10^{8}$   &$2604906 $ \\   
		$28$     &$15$   &$3,5,7$        &   $2,11, \cdots, 107$ &$ 2.444167 \times  10^{40} $   & $5.89188 \times 10^{8} $  &$6485432$ \\ 
		$27$     &$15$   &$3,5,7$        &   $2,11, \cdots, 103$ &$ 2.284268 \times  10^{38}$    & $1.24903 \times 10^{10} $ &$150567585$ \\ 
	    $26$     &$14$   &$3,5,7$        &   $2, 11, \cdots, 101$ &$2.217736 \times  10^{36}$    & $2.08433 \times 10^{10} $ &$261568064$ \\ 
        \hline\hline
		\end{tabular}
		\caption{ Number of primes with $26 \le \omega \le 35$ and number of exceptions.}
		\label{Table:tab16_4}
		\end{table}

So far, so good: however, for $5 \le \omega \le 25$ the problem rapidly becomes much harder. We repeated the above procedure and listed the results in  Table~\ref{Tab:TableFiveCons_1}. Note the increase in the level of the prime divisor tree for $\omega \leq 25$, and the rapid growth of the number of possible exceptions. It is not possible to examine each prime power on this list of exceptions because of the prohibitive computational complexity of so many cases to consider.

\begin{table}[H]
\centering
\begin{tabular}{c c c c c c}	
\hline\hline
\multicolumn{6}{c}{{\bf Values of $\omega$ and $s$ for five consecutive primitive elements }}\\ \hline 
\hline\hline
$\omega$ &$s$  &$\delta$ & \; \; \; \; \; Interval &Possible &Prime divisor\\
&    &        &                         &exceptions           &tree level\\ [0.5ex]\hline	
$34$          &$18$  & $0.00358365239643321$    &$(1.001  \times 10^{55},  3.606   \times  10^{57})$ & $1077$                  &$1$ \\ 
$33$          &$17$  & $0.03955487541801610$    &$(7.205  \times 10^{52},  2.641   \times  10^{55})$ & $1096$                  &$1$ \\ 
$32$          &$16$  & $0.07605122578297950$    &$(5.259  \times 10^{50},  6.332   \times  10^{54})$ & $36120$                 &$2$ \\ 
$31$          &$16$  & $0.01987954207276780$    &$(4.014  \times 10^{48},  9.025   \times  10^{52})$ & $337195$                &$2$ \\ 
$30$          &$15$  & $0.05924962081292530$    &$(3.161  \times 10^{46},  8.934   \times  10^{51})$ & $4.24 \times 10^{6}$  &$2$ \\ 
$29$          &$14$  & $0.10349740842354500$    &$(2.797  \times 10^{44},  2.553   \times  10^{51})$ & $1.37 \times 10^{8}$  &$2$\\ 
$28$          &$14$  & $0.04298598933316810$    &$(2.566  \times 10^{42},  1.440   \times  10^{49})$ & $5.90 \times 10^{8}$  &$3$\\ 
$27$          &$13$  & $0.08971496129578480$    &$(2.398  \times 10^{40},  2.853   \times  10^{48})$ & $1.25 \times 10^{10}$ &$3$ \\ 
$26$          &$13$  & $0.02197958084873130$    &$(2.329  \times 10^{38},  4.623   \times  10^{46})$ & $2.09 \times 10^{10}$ &$3$ \\ \hline
$25$          &$12$  & $0.07148453134378090$    &$(2.306  \times 10^{36},  3.727   \times  10^{45})$ & $5.31 \times 10^{11}$ &$4$ \\ 
$24$          &$11$  & $0.12303092309635800$    &$(2.377  \times 10^{34},  1.058   \times  10^{45})$ & $5.15 \times 10^{13}$ &$4$ \\ 
$23$          &$11$  & $0.05725947886506190$    &$(2.671  \times 10^{32},  4.749   \times  10^{42})$ & $2.06 \times 10^{13}$ &$4$ \\ 
$22$          & $ 10$ & $ 0.117500442720484$    & $  (3.217\times 10^{30} ,  9.335\times 10^{41})$  & $ 1.6754 \times 10^{14}$  &$4$\\
$21$ & $ 10$ & $ 0.0456564468258548$ & $ (4.072\times 10^{28} ,  6.002\times 10^{39})$  & $ 1.7023 \times 10^{14} $  &$4$\\
$20$ & $ 9$ & $  0.114149597510786$ & $  (5.579\times 10^{26} , 7.794\times 10^{38})$  & $ 1.0488 \times 10^{16} $  &$5$\\
$19$ & $ 9$ & $  0.0232818101414087$ & $ (7.858\times 10^{24} ,  1.814\times 10^{37})$  & $ 3.4677 \times 10^{16} $  &$5$\\
$18$ & $ 8$ & $  0.0979086758130505$ & $ (1.172\times 10^{23} ,  8.126\times 10^{35})$  & $ 1.0403 \times 10^{17} $  &$5$\\
$17$ & $ 8$ & $  0.00746209582435631$&  $(1.922\times 10^{21} ,  1.353\times 10^{35})$  & $ 1.0570 \times 10^{18} $  &$5$\\
$16$ & $ 7$ & $  0.0922078585362208$ & $ (3.258\times 10^{19} ,  6.805\times 10^{32})$  & $ 3.1356 \times 10^{18} $  &$5$\\
$15$ &  $ 6$ & $  0.186547481177730$ & $  (6.148\times 10^{17} ,  1.227\times 10^{32})$ & $ 5.0963 \times 10^{19} $   &$6$\\
$14$ & $ 6$ & $  0.0755391555533084$ & $ (1.308\times 10^{16} ,  7.201\times 10^{29})$  & $ 1.4050 \times 10^{19} $  &$6$\\
$13$ & $ 5$ & $  0.191818225320750$ & $  (3.042\times 10^{14} ,  7.814\times 10^{28})$  & $ 6.5557 \times 10^{19} $  &$6$\\ 
$12$ & $ 5$ & $  0.0506115500961032$ & $ (7.420\times 10^{12} ,  1.070\times 10^{27})$  & $ 3.6835 \times 10^{19} $  &$6$ \\
$11$ & $ 4$ & $  0.185746685231238$ & $  (2.005\times 10^{11} ,  5.136\times 10^{25})$  & $ 6.5369 \times 10^{19} $  &$6$\\
$10$ & $ 4$ & $  0.0529193607530600$&  $ (6.469\times 10^{9} ,   6.011\times 10^{23})$  & $ 2.3716 \times 10^{19} $  &$6$\\
$9$ & $ 3$ & $   0.225333153856508$ & $  (2.230\times 10^{8} ,   1.896\times 10^{22})$  & $ 2.1700 \times 10^{19} $  &$6$\\
$8$ & $ 3$ & $   0.0581090735889498$&  $ (9.699\times 10^{6} ,   2.657\times 10^{20})$  & $ 6.9939 \times 10^{18} $  &$6$\\
$7$ & $ 2$ & $   0.321266968325792$&  $  (510511 ,              4.057\times 10^{18})$   & $ 2.02852 \times 10^{18} $ &$0$\\
$6$ & $ 2$ & $   0.160839160839161$ & $  (30031 ,               1.477\times 10^{16})$   & $ 7.38642 \times 10^{15} $ &$0$\\ 
$5$ & $ 1$ & $   0.545454545454545$ & $  (2311              ,     3.831\times 10^{14})$ & $ 1.91559 \times 10^{14} $ &$0$\\ 
\hline\hline
\end{tabular}
\caption{The five consecutive primitive elements problem.}
\label{Tab:TableFiveCons_1}
\end{table}

It is possible to obtain more refined sieving inequalities by considering more congruence classes modulo small primes, along the lines of Lemma \ref{lem:main2} improving on Lemma \ref{lem:main1}. Indeed, this process enabled Cohen in \cite{Cohen1,Cohen2,Cohen3} to resolve completely the problem of two consecutive primitive elements. However, in our case, any such improvements appear marginal: a new idea is required.

\subsection*{Acknowledgements}
We are grateful to Defence Science and Technology Group (DSTG) and Australian Research Council Future Fellowship for supporting this research for the first author and for second author, respectively. In particular, we would like to thank Ralph Buchholz and Garry Hughes for their helpful discussion on the \textit{Magma} routine used. Also we would like to thank Henry Haselgrove for a useful discussion to make use of OpenMp C++ with multi-threading.

\end{document}